\documentclass[11pt,a4paper]{amsart}

\usepackage[english,francais]{babel}
\usepackage[T1]{fontenc}
\usepackage{amsmath,amssymb}
\usepackage[all]{xy}
\usepackage{graphicx}
\usepackage{fancyhdr}

\newtheorem{definition}{D\'efinition }[section]

\newtheorem{theoreme}[definition]
{Th\'eor\`eme }
\newtheorem{ex}[definition]
{Exemple }

\newtheorem{rque}[definition]{Remarque}

\newcommand{\lgw}{\longrightarrow}
\newcommand{\lgm}{\longmapsto}
\newcommand{\ovl}{\overline}

\newcommand{\Spec}{\text{Spec}\,}
\newcommand{\ord}{\text{ord}}

\newcommand{\wdh}{\widehat}

\renewcommand{\l}{\lambda}
\renewcommand{\O}{\mathcal{O}}

\newcommand{\m}{\mathfrak{m}}

\newcommand{\Z}{\mathbb{Z}}

\renewcommand{\k}{\Bbbk}

\newcommand{\R}{\mathbb{R}}

\newcommand{\K}{\mathbb{K}}

\newcommand{\N}{\mathbb{N}}

\newcommand{\Q}{\mathbb{Q}}
\renewcommand{\a}{\alpha}
\renewcommand{\b}{\beta}

\newcommand{\e}{\varepsilon}
\begin{document}
\title[Approximation diophantienne dans les corps de s\'eries]{Approximation diophantienne dans les corps de s\'eries en plusieurs variables}
\author{Guillaume Rond}

\address{Department of Mathematics\\
University of Toronto\\
Toronto, Ontario M5S 2E4\\
Canada\\
(\email{rond@picard.ups-tlse.fr})}

\maketitle

\begin{abstract}
We give here a result of diophantine approximation between $\O_N$, the ring of power series in several variables, and the completion of the valuation ring that dominates $\O_N$ for the $\m$-adic topology. We deduce from this that the Artin function of a homogenous polynomial in two variables is bounded by  an affine function, which can be interpreted in term of a {\L}ojasiewicz inequality on the wedges space.
\end{abstract}

\selectlanguage{french}

\section{Introduction}
La motivation de cet article est d'\'etudier le probl\`eme de la divisibilit\'e dans l'anneau des s\'eries formelles ou convergentes en plusieurs variables, probl\`eme qui apparait dans l'\'etude des singularit\'es d'un point de vue analytique. Si $x$ et $y$ sont deux \'el\'ements de $\k[[T]]$, nous avons toujours $x/y\in\k[[T]]$ ou $y/x\in\k[[T]]$. Cela vient du fait que $\k[[T]]$ est un anneau de valuation. La valuation consid\'er\'ee ici est la valuation $\m$-adique, o\`u $\m$ est l'id\'eal maximal de l'anneau $\k[[T]]$. Dans le cas des s\'eries en plusieurs variables, si $x$ et $y$ sont deux \'el\'ements de $\k[[T_1,...,\,T_N]]$, il est en g\'en\'eral faux que $x/y\in\k[[T_1,...,\,T_N]]$ ou $y/x\in\k[[T_1,...,\,T_N]]$. Nous pouvons donner deux exemples des cons\'equences de ce probl\`eme de ``manque de divisiblit\'e'' :\\
\\
Le premier est li\'e au probl\`eme des coins introduit par M. Lejeune-Jalabert (\cite{LJ}, voir aussi \cite{Re}), que l'on peut r\'esumer comme suit. Soit $(X,\,0)$ un germe de singularit\'e de surface et soit $(\mathfrak{X},\,0)$ sa d\'esingularisation minimale. Un arc sur $(X,\,0)$ est la donn\'ee d'un morphisme $\varphi$ de $(\Spec\k[[T]],\,0)$ vers $(X,\,0)$, et un coin sur $(X,\,0)$ est la donn\'ee d'un morphisme $\varphi$ de $(\Spec\k[[T_1,\,T_2]],\,0)$ vers $(X,\,0)$. Un arc peut toujours se relever sur un \'eclat\'e de $(X,\,0)$ de centre r\'egulier si son point g\'en\'erique n'est pas inclu dans le centre de l'\'eclatement (c'est-\`a-dire si son point g\'en\'erique peut se relever). Ceci d\'ecoule du crit\`ere de propret\'e des morphismes propres, mais nous pouvons le voir aussi ``\`a la main''. Si l'anneau des fonctions sur $(X,\,0)$ est $\k\{X_1,...,\,X_n\}/(f_1,...,\,f_p)$, un arc $\varphi$ est la donn\'ee de $n$ s\'eries $x_1(T)$,..., $x_n(T)$ qui v\'erifient $f_i(x(T))=0$ pour $1\leq i\leq p$. Si $I=(x_1,...,\,x_q)$ est le centre de l'\'eclatement et qu'il existe $j\leq q$ tel que $x_j(T)\neq 0$, nous pouvons, quitte \`a r\'eordonner les $x_j$, supposer que 
$$\ord(x_q(T)))\geq\ord(x_{q-1}(T))\geq\cdots\geq\ord(x_1(T))\neq +\infty$$
o\`u $\ord$ est la valuation $\m$-adique sur $\k[[T]]$. En l'origine de la carte $x_1\neq 0$, un relev\'e de $\varphi$ est donc donn\'e par les s\'eries 
$$x_1(T),\,x_2(T)/x_1(T),...,\,x_q(T)/x_1(T),\,x_{q+1}(T),..., \,x_n(T)$$
Le probl\`eme des coins consiste \`a savoir si un coin en position ``suffisamment g\'en\'erale'' peut toujours se relever en un coin sur  $(\mathfrak{X},\,0)$. Un coin en position "suffisamment g\'en\'erale" est un coin $\varphi(T_1,\,T_2)$ tel que l'arc $\varphi(T_1,\,0)$ soit transversal au lieu singulier de $(X,\,0)$. En g\'en\'eral un coin ne se rel\`eve pas sur un \'eclat\'e, du fait que si l'on peut ordonner les $x_j(T_1,\,T_2)$ \`a l'aide de $\ord$, la valuation $\m$-adique sur $\k[[T_1,\,T_2]]$, nous ne pouvons effectuer la division comme pr\'ec\'edemment pour les s\'eries en une variable.\\
\\
Le second exemple r\'eside dans le comportement de la fonction de Artin d'un id\'eal $I$  de $\k[[T_1,...,\,T_N]][X_1,...,\,X_n]$. Si $I=(f_1(X),...,\,f_p(X))$, cette fonction est la plus petite fonction num\'erique de $\N$ dans $\N$ qui v\'erifie la propri\'et\'e suivante  :\\
$$\forall i\in\N\ \forall x\in \k[[T_1,...,\,T_N]]^n\text{ tels que } $$
$$\left(\forall l \ \ f_l(x)\in \m^{\b(i)+1}\ \right)\,\Longrightarrow
\left(\exists \,\ovl{x}\in x+\m^{i+1} \text{ tel que } \forall l \ \ f_l(\ovl{x})=0\ \right)$$
o\`u $\m$ est l'id\'eal maximal de $\k[[T_1,...,_,T_N]]$. Son existence d\'ecoule de \cite{Ar}.\\
Si $N=1$, la fonction de Artin de $I$ est toujours born\'ee par une fonction affine \cite{Gr} et cela peut s'interpr\'eter en terme d'in\'egalit\'e de type {\L}ojasiewicz sur l'espace des arcs (cf. \cite{Hi} proposition 3.1 par exemple). La preuve de ce r\'esultat d\'ecoule du fait que nous avons une in\'egalit\'e $\b(i)\leq 2\b'(i)$ pour tout $i\in\N$, o\`u $\b$ (resp. $\b'$) est la fonction de Artin de $I$ (resp. de l'id\'eal jacobien de $I$ vu comme un id\'eal de polyn\^omes \`a coefficients dans $\k[[T_1,...,\,T_N]]$). L\`a encore, cette majoration se montre en utilisant le fait que $x$ divise $y$  dans $\k[[T]]$ si $\ord(x)\leq \ord(y)$.\\
Si $N\geq 2$, l'in\'egalit\'e pr\'ec\'edente est fausse, et la fonction de Artin d'un id\'eal $I$ comme pr\'ec\'edemment n'est pas toujours born\'ee par une fonction affine \cite{Ro}.\\
\\
Une question qui se pose alors naturellement est de savoir se qui ``se passe'' quand on divise dans $\O_N:=\Bbbk[[T_1,...,T_N]]$ (ou $\O_N:=\Bbbk\{T_1,...,T_N\}$ si $\k$ est muni d'une norme). Plus pr\'ecis\'ement nous allons comparer l'anneau $\O_N$ avec le compl\'et\'e de l'anneau de valuation qui domine $\O_N$ pour la valuation $\m$-adique, qui lui est de la forme $\K[[T]]$, avec $\K$ un  corps.
Nous allons donner ici une mani\`ere d'appr\'ehender ce probl\`eme de ``manque de divisibilit\'e'', \`a l'aide d'un th\'eor\`eme d'approximation diophantienne entre $\K_N$, le corps des s\'eries en plusieurs variables, et $\wdh{\K}_N$, son compl\'et\'e pour la topologie $\m$-adique.\\
Dans \cite{Ro}, nous avons montr\'e qu'\'etudier le comportement de la fonction de Artin des polyn\^omes homog\`enes, \`a coefficients dans l'anneau de s\'eries formelles ou convergentes en plusieurs variables sur un corps $\k$, dont le seul z\'ero est $(0,...,\,0)$, est \'equivalent \`a un tel probl\`eme d'approximation diophantienne.\\
\\
Nous montrons ici le th\'eor\`eme suivant (les notations sont donn\'ees \`a la suite) :\\
\begin{theoreme}\label{applin}
 Soit $z\in\wdh{\K}_N$ alg\'ebrique sur $\K_N$ tel que $z\notin\K_N$.   Alors il existe $a\geq 1$ et $K\geq 0$ tels que
\begin{equation}\label{dioplin}\left|z-\frac{x}{y}\right|\geq K|y|^{a},\ \forall x,\,y\in \O_N.\end{equation}\\
\end{theoreme}
Nous utilisons alors ce r\'esultat pour montrer le th\'eor\`eme suivant qui donne une r\'eponse \`a une question qui nous a \'et\'e pos\'ee par M. Hickel :
\begin{theoreme}
Soit $P(X,\,Y)$ un polyn\^ome homog\`ene en $X$ et $Y$ \`a coefficients dans $\O_N$.  Alors $P$ admet une fonction de Artin born\'ee par une fonction affine.
\end{theoreme}
\subsection{Notations }
Soient $N$ un entier positif non nul et $\k$ un corps ; nous utiliserons les notations suivantes :
\begin{itemize}
\item
$\O_N$ est l'anneau des s\'eries, en $N$ variables, formelles $\Bbbk[[T_1,...,T_N]]$ ou convergentes $\Bbbk\{T_1,...,T_N\}$ (si $\k$ est muni d'une norme) selon le cas et indiff\'eremment.
\item $\m$ est l'id\'eal maximal de $\O_N$ et  $\ord$ est la valuation $\m$-adique sur $\O_N$. Cette valuation d\'efinit une norme $|\ |$ sur $\O_N$ en posant $|x|=e^{-\ord(x)}$ et cette norme induit une topologie appel\'ee topologie $\m$-adique.
\item $V_N:=\left\{\frac{x}{y}\,/\,x,\,y\in\O_N \text{ et } \ord(x)\geq\ord(y)\right\}$, l'anneau de valuation discr\`ete qui domine $\O_N$ pour $\ord$. Nous noterons $\m'$ son id\'eal maximal.
\item  Notons $\wdh{V}_N$ le compl\'et\'e  pour la topologie $\m$-adique de $V_N$. Nous avons $\wdh{V}_N=\k\left(\frac{T_1}{T_N},...,\,\frac{T_{N-1}}{T_N}\right)[[T_N]]$. En effet cet anneau correspond au compl\'et\'e de l'\'eclatement le long de $\m$ de $\O_N$. Nous noterons $\wdh{\m}$ l'id\'eal maximal de cet anneau, $\ord$ l'extension de la valuation $\m$-adique et $|\ |$ l'extension de la norme associ\'ee. 
\item $\K_N$ et $\wdh{\K}_N$ sont respectivement les corps de fractions de $\O_N$ et de $\wdh{V}_N$. On peut remarquer que $\wdh{\K}_N$ est le compl\'et\'e de $\K$ pour la norme $|\ |$.
\end{itemize}
\begin{rque}
Il existe une th\'eorie de l'approximation diophantienne entre le corps des polyn\^omes en une variable et le corps des s\'eries en 1 variables (cf. \cite{La} pour une introduction), tr\'es similaire \`a celle entre $\Q$ et $\R$.
  La diff\'erence fondamentale entre l'approximation diophantienne que nous traitons ici et celle entre $\Q$ et $\R$ r\'eside dans le fait que les \'el\'ements non nuls de $\Z$ sont de norme sup\'erieure ou \'egale \`a 1, alors que les \'el\'ements non nuls de $\O_N$ sont tous de norme inf\'erieure ou \'egale \`a 1.  En particulier il n'existe pas de version du th\'eor\`eme de Liouville dans notre cadre (cf. remarque \ref{liouv}).
\end{rque}
\begin{rque}\label{liouv}
Dans \cite{Ro}, nous avons donn\'e une suite d'\'el\'ements $x_p\in\wdh{\K}_N$, pour $p\in\N\backslash{\{0,\,1,\,2\}}$, de degr\'e 2 sur $\K_N$, pour lesquels il existe $u_{p,\,k}$ et $v_{p,\,k}$ tels que $\left|x_p-\frac{u_{k,\,p}}{v_{k,\,p}}\right|=C_p|v_k|^{\frac{p}{2}-1}$ et $\ord(v_{p,\,k})$ tend vers $+\infty$ avec $k$, o\`u $C_p$ est une constante qui d\'epend de $p$. Nous voyons donc que, contrairement au cas des nombres r\'eels alg\'ebriques, la meilleure borne $a$ telle qu'il existe $K$ avec
$$\left|x-\frac{u}{v}\right|>K|v|^a,\ \forall u,\,v\in R$$ 
ne peut pas \^etre born\'ee par le degr\'e de l'extension de $\K_N$ par $x$. Il n'existe donc pas de version du th\'eor\`eme de Liouville pour les extensions finies de $\K_N$ dans $\wdh{\K}_N$.\\

\end{rque}
Je tiens \`a remercier ici M. Hickel et M. Spivakovsky  pour leurs conseils et commentaires \`a propos de ce travail.
\section{Preuve du th\'eor\`eme principal}
Nous nous sommes inspir\'e l\`a de la preuve de la proposition 5.1 de \cite{I2}. Nous utilisons en particulier le th\'eor\`eme d'Izumi sur les In\'egalit\'es Compl\'ementaires Lin\'eaires associ\'ees \`a l'ordre $\m$-adique sur un anneau local analytiquement irr\'eductible (cf. \cite{I2} et \cite{Re1}). Pour cela, nous allons tout d'abord introduire les notations suivantes :\\
Soit $z\in\wdh{\K}_N\backslash\K_N$ alg\'ebrique sur $\K_N$. Soit $Q(Z)$ un polyn\^ome irr\'eductible de $\O_N[Z]$ annulant $z$. Ce polyn\^ome engendre l'id\'eal des polyn\^omes de $\K_N[Z]$ s'annulant en $z$. Nous pouvons \'ecrire $Q(Z)=a_0+a_1Z+\cdots+a_dZ^d$ avec $d\geq 2$. Nous noterons $P(X,\,Y):=a_0Y^d+a_1XY^{d-1}+\cdots+a_dX^d$. Consid\'erons les extensions  suivantes :
$$\O_N\lgw \O_{N+1}\lgw \O_{N+1}/(Q(Z))=\O$$
o\`u $\O_{N+1}=\k[[T_1,...,\,T_N,\,Z]]$.
Notons $\ord$ la valuation $(T_1,...,\,T_N,\,Z)$-adique sur $\O_{N+1}$ qui \'etend la valuation $\ord$ d\'efinie pr\'ec\'edemment et $\ord_{\O}$ l'ordre  $(T_1,...,\,T_N,\,Z)$-adique sur $\O$. Consid\'erons $\ovl{Q(Z)}$, le terme initial de $Q(Z)$ pour la valuation $\ord$ sur $\O_{N+1}$. Nous avons alors 
$$Gr_{\m_{\O}}\O=\frac{\left(Gr_{\m}\O_N\otimes_{\k}\left(\k\oplus (Z)/(Z)^2\oplus(Z)^2/(Z)^3\oplus\cdots\right)\right)}{(\ovl{Q(Z)})}.$$\\
Tout d'abord, nous allons nous ramener au cas o\`u $\O$ est int\`egre et $\ovl{Q(Z)}=Z^d$. Nous serons alors sous les hypoth\`eses du th\'eor\`eme d'Izumi que nous pourrons donc utiliser.\\
Pour tout $u\in\O_N$, nous notons  $Q_u(Z)=u^da_d^{d-1}Q(Z/ua_d)$. Nous avons
$$Q_u(Z)=a_0u^da_d^{d-1}+a_1u^{d-1}a_d^{d-2}Z+\cdots+Z^d.$$
Le polyn\^ome $Q_u(Z)$ est dans $\O_N[Z]$, unitaire, et irr\'eductible comme polyn\^ome de $\K_N[Z]$ car $Q(Z)$ est irr\'eductible. Donc $Q_u(Z)$ est irr\'eductible dans $\O_N[Z]$. Fixons $u$ de telle mani\`ere \`a ce que $\ovl{Q_u(Z)}=Z^d$ et tel que $\ord(u)\geq1$ (il suffit de choisir $u$ d'ordre grand). Le polyn\^ome $Q_u(Z)$ est un polyn\^ome distingu\'e  car $\ord(u)\geq1$. Donc $Q_u(Z)$ est irr\'eductible dans $\O_{N+1}$ (\cite{To}, lemme 1.7) et $\O$ est int\`egre. Les z\'eros de $Q_u(Z)$ dans $\wdh{\K}_N$ sont les $ua_dz_i$ o\`u les $z_i$ sont les z\'eros de $Q(Z)$ dans $\wdh{\K}_N$. Soit $z$ un z\'ero de $Q(Z)$ dans $\wdh{\K}_N$. Si nous avons $\left|\frac{x}{y}-ua_dz\right|\geq K|y|^{a}$, $\forall x,\,y\in\O_N$, alors nous avons $\left|\frac{x}{y}-z\right|\geq \frac{K}{|ua_d|}|y|^{a}$, $\forall x,\,y\in\O_N$.\\
Nous pouvons donc supposer que  $\O$ est int\`egre et que $\ovl{Q(Z)}=Z^d$, ce que nous ferons \`a partir de maintenant.\\
\\
Notons $z_1,...,z_p$ les diff\'erentes solutions de $Q(Z)$ autres que $z$ dans $\wdh{\K}_N$, et fixons $x$ et $y$ dans $\O_N$. Notons $\ovl{Z}$ l'image de $Z$ dans $\O$.\\
Notons $r:=\max\{\ord(z),\,\ord(z-z_k),\, k=1,..,\,p\}$ si $p\neq 0$ et $r=\ord(z)$ sinon.
Supposons  que $\ord((x/y)-z)> r$. En particulier $\ord(x)-\ord(y)=\ord(z)$. Notons $C:=\ord(x)-\ord(y)=\ord(z)$.\\
Nous avons $Y^dQ\left(\frac{X}{Y}\right)=P(X,\,Y)$. Or $Q(\ovl{Z})=0$ dans l'anneau $\O$ qui est int\`egre, donc $Q(T)=(T-\ovl{Z})(b_{d-1}T^{d-1}+b_{d-2}T^{d-2}+\cdots+b_0)$ dans $\O[T]$, et donc $P(X,\,Y)=(X-\ovl{Z}Y)(b_{d-1}X^{d-1}+b_{d-2}X^{d-2}Y+\cdots+b_0Y^{d-1})$ dans $\O[X,\,Y]$. En d\'eveloppant l'expression pr\'ec\'edente nous voyons que 
$$\begin{array}{c}
b_{d-1}=a_d\qquad\qquad\qquad\qquad\qquad\\
b_{i}=a_{i+1}+\ovl{Z}b_{i+1}, \ 1\leq i\leq d-2\\
b_0=-a_0/\ovl{Z}\qquad\qquad\qquad\qquad\quad\end{array}$$
D'o\`u 
$$b_i:=\ovl{Z}^{d-i-1}a_d+\ovl{Z}^{d-i-2}a_{d-1}+\cdots+a_{i+1}.$$
Notons 
$$h:=b_{d-1}x^{d-1}+b_{d-2}x^{d-2}y+\cdots+b_0y^{d-1}.$$
Nous avons alors $(x-\ovl{Z}y)h=P(x,\,y)$ dans $\O$. Nous pouvons \'ecrire
$$h=\frac{P(x,\,y)-a_0y^d}{x}+\frac{P(x,\,y)-a_0y^d-a_1xy^{d-1}}{x^2}y\ovl{Z}+\cdots\qquad\qquad$$
$$\qquad\qquad+\frac{P(x,\,y)-a_0y^d-a_1xy^{d-1}-\cdots-a_{d-1}x^{d-1}y}{x^d}(y\ovl{Z})^{d-1}.$$
Notons 
$$f_i:=\frac{P(x,\,y)-a_0y^d-\cdots-a_{i}x^{i}y^{d-i}}{x^{i+1}}y^{i}$$
Nous avons alors $h=f_0+f_1\ovl{Z}+\cdots+f_{d-1}\ovl{Z}^{d-1}$ et les $f_i$ sont dans $\O_N$.\\
\\
L'anneau $\O$ \'etant int\`egre, d'apr\`es le th\'eor\`eme d'Izumi \cite{I2},  il existe deux constantes $A\geq1$ et $B\geq 0$ telles que
\begin{equation}\label{izumi}A(\ord_{\O}(x-\ovl{Z}y)+\ord_{\O}(h))+B\geq \ord_{\O}(P(x,\,y))\geq\ord(P(x,\,y)).\end{equation}
Deux cas se pr\'esentent alors :\\
\textbf{Cas 1 :}
Soit $\ord(P(x,\,y))\leq \ord(a_0y^d)$, et alors dans ce cas 
\begin{equation}\label{3}\ord(P(x,\,y))\leq d\,\ord(y)+\ord(a_0).\end{equation}
\textbf{Cas 2 :} 
Soit $\ord(P(x,\,y))> \ord(a_0y^d)$. Soit $g=g_0+g_1\ovl{Z}+\cdots+g_{d-1}\ovl{Z}^{d-1}$ avec les $g_i\in \O_N$. Alors l'image de $g$ dans $Gr_{\m_{\O}}$ est non nulle puisque $\ovl{Q(Z)}$ est de degr\'e $d$ en $Z$. Donc $\ord_{\O}(g)=\min_i\{\ord(g_i)+i\}$. En particulier, 
$$\ord_{\O}(h)\leq \ord\left(\frac{P(x,\,y)-a_0y^d}{x}\right)=d\,\ord(y)-\ord(x)+\ord(a_0).$$
D'autre part $\ord_{\O}(x-\ovl{Z}y)= \min\{\ord(x),\,\ord(y)+1\}$.
Donc, d'apr\`es l'\'equation (\ref{izumi}), nous obtenons
$$\ord(P(x,\,y))\leq Ad\,\ord(y)-A\,\ord(x)+A\,\min\{\ord(x),\,\ord(y)+1\}+A\ord(a_0)+B.$$
Nous avons alors
\begin{equation}\label{4}\ord(P(x,\,y))\leq Ad\,\ord(y)+A\ord(a_0)+B\end{equation}
\\
\textbf{Conclusion :} Si $\ord((x/y)-z)> r$, alors il existe deux constantes, $a$ et $b$ telles que 
$$\ord(P(x,\,y))\leq a\,\ord(y)+b.$$
Or $P(x,\,y)=y^dQ(x/y)$. Nous pouvons \'ecrire $Q(Z)$ sous la forme
$$Q(Z)=R(Z)\prod_{k=1}^p(Z-z_k)^{n_k}.(Z-z)^n$$
avec $R$ un polyn\^ome de degr\'e $q$ en $Z$ qui n'admet pas de z\'ero dans $\wdh{\K}_N$. Nous avons alors $\ord(R(x/y))\geq c$ pour une constante $c$ ind\'ependante de $x$ et de $y$ car $\ord(x/y)=C$ est fix\'e.
 Comme $\ord((x/y)-z)> r$, nous avons $\ord((x/y)-z_k)\leq r$. D'o\`u 
$$\ord(Q(x/y))\geq c+\sum_{k=1}^p n_k\min\{\ord(x/y),\,\ord(z_k)\}+n\,\ord((x/y)-z)$$
$$\geq c+\sum_{k=1}^p n_k\min\{C,\,\ord(z_k)\}+n\,\ord((x/y)-z).$$
Notons $D=c+\sum_{k=1}^p n_k\min\{C,\,\ord(z_k)\}$. 
Nous obtenons alors 
$$(a-d)\ord(y)+b\geq n\,\ord((x/y)-z)+D$$
o\`u encore
$$\left|\frac{x}{y}-z\right|\geq K|y|^{\a}$$
avec $K:=e^{(D-b)/n}$ et $\a:=(a-d)/n$.
Si $\ord((x/y)-z)\leq r$, nous avons alors $\left|\frac{x}{y}-z\right|\geq K'$ avec $K':=e^{-r}$. Dans tous les cas nous avons l'in\'egalit\'e voulue.$\quad\Box$\\
\\

\section{Application \`a l'\'etude de fonctions de Artin}
Nous pouvons alors donner le r\'esultat suivant :
\begin{theoreme}
Soit $P(X,\,Y)$ un polyn\^ome homog\`ene en $X$ et $Y$ \`a coefficients dans $\O_N$.  Alors $P$ admet une fonction de Artin born\'ee par une fonction affine.\\
\end{theoreme}
\begin{ex}
Soit $Q(Z)\in\O_N[Z]$ un polyn\^ome unitaire n'ayant aucune racine dans $\O_N$. Par exemple $Q(Z)=Z^d-T_1^{d+1}$ ou $Q(Z)=Z^d-(T_1^d+T_2^{d+1})$. D\'efinissons $P(X,\,Y):=Y^dQ(X/Y)$. Le polyn\^ome $P$ est homog\`ene en $X$ et $Y$ et n'admet pas d'autre solution dans $\O_N$ que $(0,\,0)$. Un tel polyn\^ome admet donc une fonction de Artin born\'ee par une fonction affine. Ceci  peut s'\'enoncer sous la forme suivante : il existe $a$ et $b$ deux constantes telles que
$$\forall x, y\in\O_N,\ \ \ord(P(x,\,y))\leq a\min\{\ord(x),\,\ord(y)\}+b.$$
En notant $\vert (x,\,y)\vert:=\max\{\vert x\vert,\,\vert y\vert\}$, ceci peut encore se r\'e\'enoncer en terme d'in\'egalit\'e de type {\L}ojasiewicz : c'est-\`a-dire qu'il existe $K>0$ et $a\geq 1$ tels que
$$\vert P(x,\,y)\vert\geq K\vert (x,\,y)\vert^a,\ \ \forall x,y\in\O_N.$$\\
La preuve du th\'eor\`eme nous permet de dire que l'on peut choisir $\a$ \'egal \`a $d$ dans le premier cas (car $Q(Z)$ n'a pas de z\'ero dans $\wdh{V}_N$). Dans le second cas $\a$ est strictement plus grand que $d$ car $Q(Z)$ admet des z\'ero dans $\wdh{V}_N$.
\end{ex}
\textbf{Preuve :}
Soit $P$ comme dans l'\'enonc\'e et $x,\,y\in \O_N$. Quitte \`a renommer les variables, nous pouvons supposer que $\ord(y)\leq\ord(x)$. Nous avons alors 
$$P(x,\,y)=y^dP\left(1,\frac{x}{y}\right) .$$
 Notons alors $Q(Y)$ le polyn\^ome $P(1,\,Y)$.\\
i) Supposons que $Q$ n'a pas de racine dans $\wdh{V}_N$. Comme $\wdh{V}_N$ est de la forme $\K[[T]]$ o\`u $\K$ est un corps et $T$ une variable formelle, d'apr\`es le th\'eor\`eme de Greenberg (cf. \cite{Gr}), $Q$ admet une fonction de Artin born\'ee par une constante $c$ et dans ce cas nous obtenons alors  $\ord\left(P\left(1,\frac{x}{y}\right)\right)< c+1$. Donc $$\ord\left(P(x,\,y)\right)<d\,\ord(y)+c+1\ .$$
ii) Si $Q$ a des racines dans $\wdh{V}_N$, toujours d'apr\`es \cite{Gr}, $Q$ admet  une fonction de Artin born\'ee par une fonction affine $i\lgm \l i+\mu$ o\`u $\l\leq d$. Remarquons que $Q$ n'a qu'un nombre fini de racines car $\wdh{V}_N$ est int\`egre. Soit $z$ un z\'ero de $Q$ dans $\wdh{\K}_N$, plus proche de $x/y$ que tous les autres z\'eros de $Q$  pour la topologie $\m$-adique. Comme $i\lgm \l i+\mu$ majore la fonction de Artin de $Q$, nous avons 
$$\ord\left(Q\left(\frac{x}{y}\right)\right)\leq \l\, \ord\left(z-\frac{x}{y}\right)+\mu.$$
$\bullet$ Si $z\in V_N$, alors $z=u/v$ avec $u$ et $v$ dans $\O_N$ premiers entre eux, et 
$$\ord\left(z-\frac{x}{y}\right)=\ord\left(\frac{uy-vx}{vy}\right).$$
Donc
$$\ord\left(P(x,\,y)\right)\leq (d-\l)\ord(y)+\l\ord(uy-vx)+(\mu-\l\ord(v)).$$
D'apr\`es le lemme d'Artin-Rees, il existe $i_0\geq 0$ ne d\'ependant que de $u$ et $v$, tel que 
\begin{equation}\label{AR}(u,\,v)\cap\m^{i+i_0}\subset (u,\,v)\m^{i}\end{equation} pour tout entier $i$ positif.
Donc si $\ord(uy-vx)\geq i+i_0$, alors il existe $\e_1$ et $\e_2$ dans $\m^{i}$ tels que $uy-vx=u\e_1-v\e_2$. En posant $\ovl{x}=x-\e_1$ et $\ovl{y}=y-\e_1$, alors $u\ovl{y}-v\ovl{x}=0$ et $\ovl{x}-x,\,\ovl{y}-y\in\m^i$. Nous choisirons dor\'enavant  une constante $i_0$ pour laquelle l'inclusion (\ref{AR}) ci-dessus est v\'erifi\'ee pour tout $(u,\,v)$, o\`u $u/v$ est une racine de $Q$ et $u$ et $v$ sont premiers entre eux. Comme $Q$ n'a qu'un nombre fini de racines, une telle constante existe.\\
$\bullet$ Si $z\notin V_N$, d'apr\`es le th\'eor\`eme \ref{applin}, il existe $a$ et $b$ tels que 
$$\ord\left(z-\frac{x}{y}\right)\leq a\,\ord{y}+b.$$
Donc nous avons
$$\ord(P(x,\,y))\leq \l\, \ord\left(z-\frac{x}{y}\right)+d\ord(y)+\mu\leq (a\l+d)\ord{y}+\l b+\mu.$$
Nous pouvons faire le m\^eme raisonnement si $\ord(y)\geq\ord(x)$. Donc dans tous les cas nous avons
$$\ord\left(P(x,\,y)\right)\leq A\min\left\{\ord(x),\,\ord(y)\right\}+B\max_{u/v\text{ z\'ero de } Q}\ord(uy-vx)+C$$
o\`u $A$, $B$ et $C$ sont des constantes, et $u/v$ est un z\'ero de $Q$ dans $\K_N$ avec $u$ et $v$ premiers entre eux dans $\O_N$.\\
Supposons que nous ayons $\ord(P(x,\,y)\geq 2\max\{A,\,B\}(i+i_0)+C$, alors nous avons soit $\min\left\{\ord(x),\,\ord(y)\right\}\geq i+i_0\geq i$, soit $\ord(uy-vx)\geq i+i_0$. Dans le premier cas, nous posons $\ovl{x}=\ovl{y}=0$. Dans le second cas, d'apr\`es le lemme d'Artin-Rees, il existe $\ovl{x}$ et $\ovl{y}$ tels que $u\ovl{y}-v\ovl{x}=0$ et $\ovl{x}-x,\,\ovl{y}-y\in\m^i$. Dans tous les cas $P(\ovl{x},\ovl{y})=0$ et $\ovl{x}-x,\,\ovl{y}-y\in\m^i$. C'est-\`a-dire que $P$ admet une fonction de Artin born\'ee par la fonction affine $i\lgm 2\max\{A,\,B\}(i+i_0)+C$.$\quad\Box$\\


\begin{thebibliography}{}
\bibitem[Ar]{Ar} M. Artin, Algebraic approximation of structures over complete local rings, \textit{Publ. Math. IHES}, \textbf{36}, (1969), 23-58.
\bibitem[Gr]{Gr} M. J. Greenberg, Rational points in henselian discrete valuation rings, \textit{Publ. Math. IHES}, \textbf{31}, (1966), 59-64.
\bibitem[Hi]{Hi} M. Hickel, Calcul de la fonction d'Artin-Greenberg d'une branche plane, \textit{Pacific J. of Math.}, \textbf{213}, (2004), 37-47.
\bibitem[Iz]{I2} S. Izumi, A measure of integrity for local analytic algebras, \textit{Publ. RIMS, Kyoto Univ.}, \textbf{21}, (1985), 719-736.
\bibitem[La]{La} A. Lasjaunias,  A survey of Diophantine approximation in fields of power series, \textit{Monatsh. Math.}, \textbf{130},  (2000),  no. 3, 211-229.
\bibitem[LJ]{LJ} M. Lejeune-Jalabert, Arcs analytiques et r\'esolution minimales des singularit\'es des surfaces quasihomog\`enes, \textit{Lectures Notes in Math.}, \textbf{777}, (1980), 303-336.
\bibitem[Na]{Na} J. Nash, Arcs structure of singularities, \textit{Duke Math. J.}, \textbf{81}, (1995), 31-38.
\bibitem[Ree]{Re1} D. Rees, Izumi's theorem, Commutative algebra (Berkeley, CA, 1987),  407--416, \textit{Math. Sci. Res. Inst. Publ.}, \textbf{15}, Springer, New York, (1989).
\bibitem[Reg]{Re} A. Reguera,  Image of the Nash map in terms of wedges, \textit{C. R. Math. Acad. Sci. Paris}, \textbf{338},  (2004),  no. 5, 385-390.
\bibitem[Ro]{Ro} G. Rond, Contre-exemple \`a la lin\'earit\'e de la fonction de Artin, pr\'epublication ArXiv, (2004).
\bibitem[To]{To} J.-C. Tougeron, \textit{Id\'eaux de fonctions diff\'erentiables}, Ergebnisse der Mathematik und ihrer Grenzgebiete, Band 71. Springer-Verlag, Berlin-New York, (1972).
\bibitem[Wa]{Wa} J. J. Wavrik, A theorem on solutions of analytic equations with applications to deformations of complex structures, \textit{Math. Ann.}, \textbf{216}, (1975), 127-142.
\end{thebibliography}
\end{document}